\documentclass[12pt,reqno,final]{amsart}
\usepackage{cirrhosis}
\hypersetup{pdfpagemode=UseOutlines,colorlinks=false,pdfpagelayout=SinglePage,pdfstartview=FitH,bookmarksopen=true}

\begin{document}

\title{Holomorphic Koszul-Brylinski Homology}
\author[M.~Sti\'enon]{Mathieu Sti\'enon}
\address{Universit\'e Paris Diderot, Institut de Math\'ematiques de Jussieu (UMR~CNRS~7586), B\^atiment Chevaleret, Case~7012, 
75205~Paris~cedex~13, France}
\email{stienon@math.jussieu.fr}
\address{Pennsylvania State University, Department of Mathematics, 109 McAllister Building, University Park, PA 16802, United States}
\email{stienon@math.psu.edu}
\begin{abstract}
In this note, we study the Koszul-Brylinski homology of \emph{holomorphic} Poisson manifolds. 
We show that it is isomorphic to the cohomology of a certain \emph{smooth} complex Lie algebroid with values in the Evens-Lu-Weinstein duality module. 
As a consequence, we prove that the Evens-Lu-Weinstein pairing on Koszul-Brylinski homology is nondegenerate.
Finally we compute the Koszul-Brylinski homology for Poisson structures on $\cposq$.
\end{abstract}
\maketitle
\section{Introduction}

In~\cite{Brylinski}, Brylinski introduced a homology theory for Poisson manifolds, 
which is nowadays called Koszul-Brylinski homology. 
Evens, Lu \& Weinstein~\cite{ELW} and Xu~\cite{Xu} proved independently 
that, for unimodular Poisson manifolds, the Koszul-Brylinski homology is 
(up to a change of degree) isomorphic to the Lichnerowicz-Poisson cohomology~\cite{Lichnerowicz}. 
And Evens, Lu \& Weinstein introduced a pairing on Koszul-Brylinski homology groups. 
In this note, we study the Koszul-Brylinski homology of \emph{holomorphic} Poisson manifolds. 
Koszul-Brylinski homology is defined as the hypercohomology of the complex of sheaves
\[ \cdots\xto{\delpi}\hOO^{i+1}\xto{\delpi}\hOO^i\xto{\delpi}
\hOO^{i-1}\xto{\delpi}\cdots ,\]
where $\delpi=i_{\pi}\rond\partial-\partial\rond i_{\pi}$.
As is explained in~\cite{LSX}, any holomorphic Poisson manifold gives rise 
to a holomorphic Lie algebroid structure $\txspi$ on the holomorphic vector bundle $\txs$,
which in turn induces a complex Lie algebroid structure $\tzu\bowtie\tuzspi$ 
on the complex vector bundle $\tzu\oplus\tuzs$.
We show that the cohomology of this complex Lie algebroid with values in the Evens-Lu-Weinstein duality module is isomorphic to the Koszul-Brylinski homology. 
As a consequence, we prove that the Evens-Lu-Weinstein pairing on Koszul-Brylinski homology is nondegenerate.
We also introduce the Euler characteristic for the Koszul-Brylinski homology of a Poisson manifold
and show that it coincides with the signed Euler characteristic of the manifold. 
Finally we compute the Koszul-Brylinski homology for Poisson structures on $\cposq$.
We refer the reader to the works of Etingof \& Ginzburg~\cite{EG} and Pichereau~\cite{Pichereau} for more on the Koszul-Brylinski homology of algebraic Poisson varieties.

\section{Holomorphic Lie algebroid cohomology}

Let $A$ be a holomorphic Lie algebroid over a complex manifold $X$: i.e. $A\to X$ is a holomorphic vector bundle whose sheaf of holomorphic sections $\shs{A}$ is endowed with a Lie bracket $\lb{\cdot}{\cdot}:\shs{A}\times\shs{A}\to\shs{A}$, and there exists a holomorphic bundle map $A\xto{a}\tx$, called anchor, which induces a morphism of sheaves of $\hFF$-modules
 $\shs{A}\xto{a}\hXX$ such that 
\begin{gather} a(\lb{s_1}{s_2})=\lb{a(s_1)}{a(s_2)}, \qquad\forall s_1,s_2\in\shs{A} ; \\ 
\lb{s_1}{f s_2}=\big(a(s_1)f\big) s_2+f\lb{s_1}{s_2}, \qquad\forall s_1,s_2\in\shs{A},\;f\in\hFF .\end{gather} 

This holomorphic Lie algebroid structure gives rise to a  complex of sheaves: 
\[ \cdots\xto{d_A}\Omega_A^{k-1}\xto{d_A}
\Omega_A^{k}\xto{d_A}
\Omega_A^{k+1}\xto{d_A}\cdots ,\] 
where $\Omega_A^k$ stands for the sheaf of holomorphic sections of
 the holomorphic vector bundle $\wedge^k A^*$,
 and $d_A$ is given by the usual Cartan formula.
By definition~\cites{ELW, LSX}, the holomorphic Lie algebroid cohomology of $A$ (with trivial coefficients) is the hypercohomology of this complex of sheaves: \[ H^*(A,\CC) := \HH^*(X,\shss{A}^\bullet) .\]

A holomorphic vector bundle $E\to X$ (with sheaf of holomorphic functions $\shs{E}$)
is said to be a module over the holomorphic Lie algebroid $A$,
if there is a morphism of sheaves (of $\CC$-modules)
\[ \shs{A}\otimes\shs{E}\to\shs{E}:V\otimes s\mapsto\nabla_V s \]
such that, for any open subset $U\subset X$, the relations
\begin{gather*}
\nabla_{fV} s=f\nabla_V s \\
\nabla_V (fs)=\big(\rho(V)f\big) s +f \nabla_V s \\
\nabla_V\nabla_W s-\nabla_W\nabla_V s=\nabla_{\lie{V}{W}} s
\end{gather*}
are satisfied $\forall f\in\hFF(U)$, $\forall V,W\in\shs{A}(U)$
and $\forall s\in\shs{E}(U)$. Such a morphism $\nabla$ is called a representation of $A$ on $E$.
Given an $A$-module $E\to X$,
 one can form the complex of sheaves 
\beq{Alabama} \cdots\xto{d_A^\nabla}\shss{A}^{k-1}
\otimes_{\hFF}\shs{E}\xto{d_A^\nabla} \shss{A}^{k}\otimes_{\hFF}\shs{E}\xto{d_A^\nabla}
\shss{A}^{k+1}\otimes_{\hFF}\shs{E}
\xto{d_A^\nabla}\cdots .\eeq
By definition, the Lie algebroid cohomology of $A$ with values in $E$ is the hypercohomology of this complex of sheaves: 
\[ H^*(A,E):=\HH^*(X,\shss{A}^\bullet\otimes_{\hFF}\shs{E}) .\]

Given a holomorphic Lie algebroid $A$ with anchor $a$, we define $a^{1,0}=\frac{1-i J}{2}\smalcirc a: A\to T_\CC X$. 
Here $J$ stands for the almost complex structure $J:\tx\to\tx$ of the complex manifold $X$.
Of course, for any holomorphic function $f\in\hFF(U)$, we have 
$a^{1,0}(V)f=a(V)f$, for all $V\in\sections{U,A}$. 
Now regard $A$ as a complex vector bundle.
The Lie bracket, which was defined so far only on the sheaf of holomorphic sections of $A$, extends naturally to
all smooth sections through the Leibniz rule 
\[ \lb{s_1}{f s_2}=\big(a^{1,0}(s_1)f\big) s_2+f\lb{s_1}{s_2}, \qquad\forall s_1,s_2\in\sections{A},\;f\in\cinf{X,\CC} ,\]
with $a^{1,0}$ substituted to $a$. 
We use the symbol $\auz$ to denote the resulting complex Lie algebroid structure on $A$~\cite{LSX}.

Now recall that the complex vector bundle $\tzu$ is endowed with a canonical complex Lie algebroid structure whose Lie bracket is completely determined by the relation $\lie{\partial_{\cc{z_j}}}{\partial_{\cc{z_k}}}=0$ and the anchor, which is simply the injection 
$\tzu\injection\tcx$.

\begin{prop}[\cite{LSX}*{Theorems~4.2 and~4.8}]
\label{Indianapolis}
If $A$ is a holomorphic vector bundle with anchor $a$ over a complex manifold $X$, 
there exists a unique complex Lie algebroid structure on the complex vector bundle 
$T^{0, 1}_X\oplus A^{1, 0}$ with anchor $\abt(X^{0,1}+\xi)=X^{0,1}+a^{1,0}(\xi)$
such that $\lb{\ahXX}{\shs{A}}=0$ and both $T^{0,1}_X$ and $A^{1, 0}$ are Lie subalgebroids.
\end{prop}

This complex Lie algebroid is denoted $\tzu\bowtie\auz$. 
The pair $(\tzu,\auz)$ is an example of matched pair~\cites{Lu, Mokri, LSX}. 

\begin{thm}[\cite{LSX}*{Lemma~4.16 and Theorem~4.19}]
\label{Des_Moines}
Let $A\to X$ be a holomorphic Lie algebroid and $E\to X$ a complex vector bundle. 
Then $E$ is a module over the holomorphic Lie algebroid $A$ if, and only if,
$E$ is a module over the complex Lie algebroid $\tzu\bowtie \auz$. 
Moreover, we have \[ H^*(A, E) \isomorphism H^*(\tzu\bowtie \auz, E) .\]
\end{thm}

Note that the complex Lie algebroid $\tzu\bowtie\auz$
is an elliptic Lie algebroid in the sense
of Block~\cite{Block}. That is, $\Re\rond\abt$ is surjective. 
Therefore, when $X$ is compact, the cohomology groups $H^*(\tzu\bowtie \auz, E)$ are finite dimensional 
and we can consider the Euler characteristic 
\begin{equation}
\Euler{A,E}=\sum_i (-1)^i \dim H^i(A,E)
.\end{equation}

\begin{prop}
\label{Baton_Rouge}
Let $A\to X$ be a holomorphic Lie algebroid and $E$ an $A$-module.
Assume that $X$ is compact. Then
\[ \Euler{A,E}=\sum_i (-1)^i \Euler{X,\wedge^iA^*\otimes E} ,\]
where $\Euler{X,\wedge^iA^*\otimes E}$ is the Euler characteristic 
 of the holomorphic bundle $\wedge^iA^*\otimes E$. 
\end{prop}

\begin{proof}
By definition, $H^*(A,E)$ is isomorphic to the hypercohomology 
$\HH^*(X,\OO^{\bullet}_A\otimes_{\hFF}\shs{E})$ 
of the complex of sheaves \eqref{Alabama}, which, according to 
Theorem~\ref{Des_Moines}, is computed by the total cohomology $H^n(\tzu\bowtie\auz,E)$ of the double complex \begin{equation*}
\xymatrix{
\vdots & \vdots & \vdots & \\
\OA{0}{2} \ar[u]^{d_{\auz}^\nabla} \ar[r]^\delbar 
& \OA{1}{2} \ar[u]^{d_{\auz}^\nabla} \ar[r]^\delbar
& \OA{2}{2} \ar[u]^{d_{\auz}^\nabla} \ar[r]^\delbar 
& \cdots \\
\OA{0}{1} \ar[u]^{d_{\auz}^\nabla} \ar[r]^\delbar 
& \OA{1}{1} \ar[u]^{d_{\auz}^\nabla} \ar[r]^\delbar
& \OA{2}{1} \ar[u]^{d_{\auz}^\nabla} \ar[r]^\delbar 
& \cdots \\
\OA{0}{0} \ar[u]^{d_{\auz}^\nabla} \ar[r]^\delbar 
& \OA{1}{0} \ar[u]^{d_{\auz}^\nabla} \ar[r]^\delbar
& \OA{2}{0} \ar[u]^{d_{\auz}^\nabla} \ar[r]^\delbar 
& \cdots }
\end{equation*}
where $\hOO^{i,j}=\sections{\wedge^i\tuz\otimes\wedge^j\tzu}$ 
and $\ssof{A}^{k,l}=\sections{\wedge^k(\auz)^*\otimes\wedge^l(\azu)^*\otimes E}$. 

Set $C^{p,q}=\OA{p}{q}$ and $C^n=\bigoplus_{p+q=n}C^{p,q}$.  
The spectral sequence induced by the filtration 
$F_q(C^n)=\bigoplus_{\substack{\tilde{q}\geq q \\ \tilde{p}+\tilde{q}=n}} C^{\tilde{p},\tilde{q}}$ 
of $C^{\bullet}$ starts with $E^{p,q}_0=C^{p,q}$, $d^{p,q}_0=\delbar$ and 
$E^{p,q}_1=H^p(C^{\bullet,q},\delbar)$, and converges to 
$H^n(\tzu\bowtie\auz,E)$. 

Since the Euler characteristic of $E^{p,q}_r$ does not change from one sheet to the next, we have 
\begin{align*}
\Euler{A,E}&=\sum_n(-1)^n\dim H^n(A,E) \\ 
&=\sum_n(-1)^n\dim H^n(\tzu\bowtie\auz,E) \\ 
&=\sum_n(-1)^n\dim \big(\bigoplus_{p+q=n}E^{p,q}_\infty\big) \\ 
&=\sum_n(-1)^n\dim \big(\bigoplus_{p+q=n}E^{p,q}_1\big) \\ 
&=\sum_n(-1)^n\dim \big(\bigoplus_{p+q=n} H^p(C^{\bullet,q},\delbar)\big) \\ 
&=\sum_q(-1)^q \big(\sum_p(-1)^p\dim H^p(C^{\bullet,q},\delbar)\big) \\ 
&=\sum_q(-1)^q\Euler{X,\wedge^q A^*\otimes E} 
. \qedhere \end{align*}
\end{proof}
 
\section{Holomorphic Poisson manifolds}

A holomorphic Poisson manifold is a complex manifold $X$ whose sheaf of holomorphic funcions $\hFF$ is a sheaf of Poisson algebras. 
By a sheaf of Poisson algebras over $X$, we mean that, for each open subset $U\subset X$, the ring $\hFF(U)$ is endowed with a Poisson bracket such that all restriction maps $\hFF(U)\to\hFF(V)$ (for arbitrary open subsets $V\subset V\subset X$) are morphisms of Poisson algebras. Moreover, given an open subset $U\subset X$, an open covering $\{U_i\}_{i\in I}$ of $U$, and a pair of functions $f,g\in\hFF(U)$, the local data $\pb{f|_{U_i}}{g|_{U_i}}$ ($i\in I$) glue up and give $\pb{f|_U}{g|_U)}$ if they coincide on the overlaps $U_i\cap U_j$. On a given complex manifold $X$, the holomorphic Poisson structures are in one-to-one correspondence with the sections $\pi\in\sections{\wedge^2\tuz}$ such that $\cc{\partial}\pi=0$ and $\schouten{\pi}{\pi}=0$. The Poisson bracket on functions and the bivector field are related by the formula $\pi(\partial f,\partial g)=\pb{f}{g}$, where $f,g\in\hFF$.

Given a holomorphic Poisson bracket 
\[ \hFF\otimes_{\CC}\hFF\to\hFF:(f,g)\mapsto\pb{f}{g} ,\] 
the formula 
\beq{Kansas} \lb{f_1\;dg_1}{f_2\;dg_2}=f_1 X_{g_1}(f_2)\;dg_2-f_2 X_{g_2}(f_1)\;dg_1
+f_1 f_2 \; d\pb{g_1}{g_2} ,\eeq
where $f_1,f_2,g_1,g_2\in\hFF$, defines a Lie bracket on $\hOO$. 
Here $\ham{f}\in\hXX$ denotes the derivation \[ \ham{f}:\hFF\to\hFF:g\mapsto\pb{f}{g} \] of $\hFF$ associated to the holomorphic function $f\in\hFF$. 
Since $\sections{\tuzs}=\FF\hOO$, the bracket on $\hOO$ extends to $\sections{\tuzs}$ by the Leibniz rule: 
\[ \lb{f dz_k}{g dz_l}=f\ham{z_k}(g)dz_l-g\ham{z_l}(f)dz_k+fg\, d\pb{z_k}{z_l} ,\] for all $f,g\in\FF$. 
If the bivector field associated to the Poisson bracket on $\hFF$ is $\pi\in\hXX^2\subset\sections{\wedge^2\tuz}$, then 
the Lie bracket is given by 
\begin{equation*}
\lb{\alpha}{\beta}=L_{\pi\diese\alpha}\beta
-L_{\pi\diese\beta}\alpha-\partial(\pi(\alpha,\beta)) 
,\quad\forall\alpha,\beta\in\sections{\tuzs} 
.\end{equation*}
Once its sheaf of sections $\hOO$ has been endowed with this Lie bracket, the cotangent bundle $\txs$ becomes a holomorphic Lie algebroid with anchor map $\pi^\sharp:\txs\to\tx$, which we refer to by the symbol $\txspi$.
By Proposition~\ref{Indianapolis}, we can associate to it the complex Lie algebroid $\tzu\bowtie\tuzspi$.

The complex Lie algebroid structure on $\tzu\bowtie\tuzspi$ is characterized as follows: 
the anchor is $\id_{\tzu}\oplus\pi\diese:\tzu\oplus\tuzs\to\tcx$, 
and the Lie bracket on $\Gamma (\tzu\bowtie\tuzspi)$ satisfies $\lb{\ahXX}{\hOO}=0$, coincides with the Lie bracket of vector fields on $\ahXX$ and with the bracket defined by~\eqref{Kansas} on $\hOO$~\cite{LSX}.

\section{Holomorphic Koszul-Brylinski homology}

Let $\hXX^k$ and $\hOO^k$ denote the sheaves of holomorphic sections of
 $\wedge^k \tx$ and $\wedge^k\txs$, respectively. 

The Koszul-Brylinski operator 
$\delpi:\hOO^k\to\hOO^{k-1}$ is defined as 
$\delpi:=\ii{\pi}\del-\del\ii{\pi}$, where 
$\del:\hOO^k\to\hOO^{k+1}$ is the holomorphic exterior differential (i.e. the Dolbeault operator) and 
$\ii{\pi}:\hOO^k\to\hOO^{k-2}$ is the contraction with the holomorphic Poisson bivector field $\pi$ 
\cites{Brylinski,Koszul}. 
The operator $\delpi$ satisfies $\delpi^2=0$, $\delpi d+d\delpi =0$, and 
\[ \delpi(\alpha\wedge\beta)= \delpi \alpha \wedge \beta + (-1)^k \alpha \wedge \delpi \beta+  (-1)^k [\alpha, \beta],  
\quad \forall \alpha\in\hOO^k, \beta\in\hOO^l .\]

\begin{defn}
Let $(X,\pi)$ be a holomorphic Poisson manifold.
Its Koszul-Brylinski homology is the hypercohomology of the complex of sheaves
\begin{equation} \label{Rhode_Island} 
\cdots\xto{\delpi}
\hOO^{k+1}\xto{\delpi}\hOO^k\xto{\delpi}\hOO^{k-1}
\xto{\delpi}\cdots \end{equation}
which is denoted $H_*(X,\pi)$.
\end{defn}

\begin{rmk}
If $\pi=0$, we have $H_k(X,\pi)\isomorphism\bigoplus_{j-i=n-k} H^j(X,\hOO^i)$.
\end{rmk}

As was pointed out earlier, a holomorphic Poisson manifold $(X,\pi)$ automatically gives rise 
to a holomorphic Lie algebroid structure $\txspi$. 
The Lichnerowicz-Poisson cohomology $H^*(X,\pi;E)$ of $(X,\pi)$ with coefficients in a $\txspi$-module $E$ is defined to be
the Lie algebroid cohomology of $\txspi$ with coefficients in the module $E$, 
i.e. the hypercohomology of the complex of sheaves 
\begin{equation*}
\cdots \xto{d_{\pi}^{\nabla}}
\hXX^{k-1}\otimes_{\hFF}\shs{E} \xto{d_{\pi}^{\nabla}}
\hXX^{k}\otimes_{\hFF}\shs{E} \xto{d_{\pi}^{\nabla}}
\hXX^{k+1}\otimes_{\hFF}\shs{E} \xto{d_{\pi}^{\nabla}}
\cdots 
.\end{equation*}

In particular, when $E$ is the trivial module $X\times\CC\to X$,
the associated differential complex is 
\[ \cdots\xto{d_\pi}\hXX^{k-1}\xto{d_\pi}
\hXX^{k}\xto{d_\pi}\hXX^{k+1}\xto{d_\pi}\cdots .\]
One has $d_\pi V=\lb{\pi}{V}$. 
The hypercohomology of this complex of sheaves is the holomorphic Lichnerowicz-Poisson cohomology 
$H^*(X,\pi)$ of the holomorphic Poisson manifold $(X,\pi)$~\cite{LSX}. 

Assuming $X$ compact, let 
\begin{equation}
\EulerLP{X,\pi;E}=\sum_i (-1)^i \dim H^i(X,\pi;E)
\end{equation}
be the Euler characteristic of the Lichnerowicz-Poisson cohomology $H^*(X,\pi;E)$.

\begin{prop}\label{Kalamazoo}
If $(X,\pi)$ is a compact holomorphic Poisson manifold, then 
\[ \EulerLP{X,\pi;E}=\sum_i (-1)^i \Euler{X,\wedge^i\tx\otimes E} ,\]
where $\Euler{X,\wedge^i\tx\otimes E}$ stands for the usual Euler characteristic of the holomorphic bundle $\wedge^i\tx\otimes E$.
\end{prop} 

\begin{proof}
By definition, we have $H^k(X,\pi;E)=H^k(\txspi,\wedge^n\txs)$, whence 
\[ \EulerLP{X,\pi;E}=\Euler{\txspi,\wedge^n\txs} .\] 
Therefore, it suffices to apply Proposition~\ref{Baton_Rouge} to the Lie algebroid $A=\txspi$ and its module $E=\wedge^n\txs$ to conclude.
\end{proof}

A result of Evens, Lu \& Weinstein (transposed to the holomorphic setting)
 asserts that, if $A\to X$ is a holomorphic Lie algebroid
 with $\dim_{\CC}X=n$ and $\rk_{\CC}A=r$, the holomorphic vector bundle 
$Q_A=\wedge^r A\otimes\wedge^n\txs$ is naturally a module over $A$. 
When the holomorphic Lie algebroid $A$ is the cotangent 
bundle $\txspi$ of a holomorphic Poisson manifold $(X,\pi)$,
 we have $Q_A=\wedge^n \txs\otimes\wedge^n\txs$.
 Its square root $\sqrt{Q_A}=\wedge^n\txs$ is also an $A$-module; the representation is the map 
\[ \hOO\otimes\hOO^n\to\hOO^n:
\alpha\otimes\omega\mapsto\nabla_\alpha\omega \] 
such that $\nabla_{df}\omega=\ld{\ham{f}}\omega$, for all $f\in\hFF$
 and $\omega\in\hOO^n$.
Here $\hOO$ and $\hOO^n$ are the sheaves of holomorphic sections of
 $\txs$ and $\wedge^n\txs$ respectively. 
Hence, we obtain the complex of sheaves 
\begin{equation} \cdots \xto{d_{\pi}^{\nabla}}
\hXX^{k-1}\otimes_{\hFF}\hOO^n \xto{d_{\pi}^{\nabla}} 
\hXX^{k}\otimes_{\hFF}\hOO^n \xto{d_{\pi}^{\nabla}} 
\hXX^{k+1}\otimes_{\hFF}\hOO^n \xto{d_{\pi}^{\nabla}} 
\cdots .\label{New_Jersey}\end{equation}

An argument of Evens, Lu \& Weinstein (see~\cite{ELW}*{Equation~(22)}) 
adapted to the holomorphic context shows that the isomorphism of sheaves of $\hFF$-modules 
\[ \tau:\hXX^k\otimes_{\hFF}\hOO^n\to\hOO^{n-k}:
X\otimes\alpha\mapsto\ii{X}\alpha \] 
is in fact an isomorphism between the complexes of sheaves~\eqref{New_Jersey} and~\eqref{Rhode_Island}: 
\beq{Iowa} \xymatrix{ 
\hOO^n \ar[r] \ar[d]_{\id} & 
\cdots \ar[r] & \hXX^k\otimes_{\hFF}\hOO^n \ar[r]^{d_{\pi}^{\nabla}} \ar[d]_{\tau} & 
\hXX^{k+1}\otimes_{\hFF}\hOO^n \ar[r] \ar[d]_{\tau} & 
\cdots \ar[r] & \hXX^n\otimes_{\hFF}\hOO^n \ar[d]_{\tau} \\ 
\hOO^n \ar[r] & 
\cdots \ar[r] & \hOO^{n-k} \ar[r]_{(-1)^{k+1}\delpi} & 
\hOO^{n-k-1} \ar[r] & 
\cdots \ar[r] & \hOO^0 
} \eeq

This isomorphism of complexes of sheaves induces an isomorphism of the corresponding sheaf cohomologies.
Thus we obtain the following theorem, which is a holomorphic analogue of a result of 
Evens, Lu \& Weinstein~\cite{ELW}*{Corollary~4.6}. 

\begin{thm}\label{Vermont}
For any holomorphic Poisson manifold $(X,\pi)$, the chain map $\tau$ induces an isomorphism 
\[ H^k\big(X,\pi;\wedge^n\txs\big) \xto{\isomorphism} H_{2n-k}(X,\pi) .\]
\end{thm}

Assume that $(X,\pi)$ is a compact holomorphic Poisson manifold. Let 
\begin{equation}
\EulerKB{X,\pi}=\sum_i (-1)^i \dim H_i(X,\pi)
\end{equation}
be the Euler characteristic of the Koszul-Brylinski homology.

\begin{thm}\label{Silver_Springs}
For a compact holomorphic Poisson manifold $(X,\pi)$, we have 
\[ \EulerKB{X,\pi}=(-1)^n\Euler{X} ,\] 
where $\Euler{X}$ denotes the standard Euler characteristic of $X$.
\end{thm}

\begin{proof}
We have 
\begin{align*} 
\EulerKB{X,\pi} &= \EulerLP{X,\pi;\wedge^n\txs} &&\text{by Theorem~\ref{Vermont}} \\ 
&= \sum_i(-1)^i\Euler{X,\wedge^i\tx\otimes\wedge^n\txs} &&\text{by Proposition~\ref{Kalamazoo}} \\ 
&= (-1)^n\sum_j(-1)^j\Euler{X,\wedge^j\txs} && \\ 
&= (-1)^n\Euler{\tx,\CC} &&\text{by Proposition~\ref{Baton_Rouge}} 
.\end{align*}
Of course, since $\hOO^{\bullet}\xto{\partial}\hOO^{\bullet+1}$ and 
$\sections{\wedge^{\bullet}\tcxs}\xto{d}\sections{\wedge^{\bullet+1}\tcxs}$ are two acyclic resolutions of the locally constant sheaf $\CC$ over $X$, we have 
\begin{align*} 
\Euler{\tx,\CC} &= \sum_i(-1)^i\dim H^i\big(\sections{\hOO^{\bullet}}\xto{\partial}\sections{\hOO^{\bullet+1}}\big) \\ 
&= \sum_i(-1)^i\dim H^i\big(\sections{\wedge^{\bullet}\tcxs}\xto{d}\sections{\wedge^{\bullet+1}\tcxs}\big) =\Euler{X}  
.\qedhere\end{align*} 
\end{proof}

\begin{defn}
A holomorphic Poisson manifold $(X,\pi)$ is said to be unimodular
if $\wedge^n\txs$ is isomorphic, as a $\txspi$-module, to the trivial module $\CC$.
\end{defn}

The notion of modular class was introduced independently by Brylinski \& Zuckerman~\cite{BZ} for holomorphic Poisson manifolds, and by Weinstein~\cite{Weinstein} for real Poisson manifolds.
For the relation between Calabi-Yau algebras and unimodular Poisson structures, see~\cite{Dolgushev}.

From the definition, it is clear that a holomorphic Poisson manifold 
$(X,\pi)$ is unimodular if and only if there exists a global
holomorphic section $\omega \in \hOO^n$ such that
the vector field $H\in\hXX$ defined by
\[ \nabla_{df}\omega=L_{X_f}\omega=H(f)\cdot\omega \qquad (f\in\hFF) \] 
is a holomorphic Hamiltonian vector field. 

\begin{prop}
For a unimodular holomorphic Poisson manifold $(X,\pi)$, the chain map $\tau$ induces an isomorphism  
\[ H^k(X,\pi)\xto{\isomorphism} H_{2n-k}(X,\pi) .\]
\end{prop} 

\section{Koszul-Brylinski double complex}

In this section, we describe a double complex computing the Koszul-Brylinski homomology.

\begin{thm}
The Koszul-Brylinski homology of a holomorphic
Poisson manifold $(X,\pi)$ is isomorphic to the
total cohomology of the double complex
\begin{equation*}
\xymatrix{
\cdots \ar[r] & \hOO^{n-k+1,0} \ar[d]_{\delbar} \ar[r]^{(-1)^{k}\delpi} &
\hOO^{n-k,0} \ar[d]_{\delbar} \ar[r]^{(-1)^{k+1}\delpi} &
\hOO^{n-k-1,0} \ar[d]_{\delbar} \ar[r] & \cdots \\
\cdots \ar[r] & \hOO^{n-k+1,1} \ar[d]_{\delbar} \ar[r]^{(-1)^{1+k}\delpi} &
\hOO^{n-k,1} \ar[d]_{\delbar} \ar[r]^{(-1)^{1+k+1}\delpi} &
\hOO^{n-k-1,1} \ar[d]_{\delbar} \ar[r] & \cdots \\
\cdots \ar[r] & \hOO^{n-k+1,2} \ar[d]_{\delbar} \ar[r]^{(-1)^{2+k}\delpi} &
\hOO^{n-k,2} \ar[d]_{\delbar} \ar[r]^{(-1)^{2+k+1}\delpi} &
\hOO^{n-k-1,2} \ar[d]_{\delbar} \ar[r] & \cdots \\
& \vdots & \vdots & \vdots
} \end{equation*}
\end{thm}

\begin{proof}
According to Theorem~\ref{Des_Moines}, we have
\begin{equation*}
H_*(X,\pi) \isomorphism
H^*(\tzu\bowtie\tuzspi,\wedge^n\tuzs) 
.\end{equation*}
The r.h.s.\ is the Lie algebroid cohomology of
$\tzu\bowtie\tuzspi$ with coefficients in
the module $\wedge^n\tuzs$.
Moreover, the representation of the complex Lie algebroid 
$\tzu\bowtie\tuzspi$ on $\wedge^n\tuzs$ is the map
\[ \sections{\tzu\oplus\tuzs}\otimes\sections{\wedge^n\tuzs} 
\to\sections{\wedge^n\tuzs}:(X+\xi,\omega)\mapsto\nabla_{X+\xi}\omega \] 
defined by 
\begin{gather*}
\nabla_{\pzb{k}}(f\;dz_1\wedge\cdots\wedge dz_n) = 
\tfrac{\partial f}{\partial \cc{z}_k} \;dz_1\wedge\cdots\wedge dz_n 
\\ 
\nabla_{dz_l}(f\;dz_1\wedge\cdots\wedge dz_n) =
\ld{\ham{z_l}}(f\;dz_1\wedge\cdots\wedge dz_n) 
\end{gather*}
(for all $f\in\FF$). 

Consider the complex 
\beq{Virginia} 
\Gamma\Big(\wedge^m\big(\tzu\oplus\tuzs\big)^*\otimes\wedge^n\tuzs\Big) 
\xto{\dee^{\nabla}_{\bowtie}} 
\Gamma\Big(\wedge^{m+1}\big(\tzu\oplus\tuzs\big)^*\otimes\wedge^n\tuzs\Big) 
\eeq
Set $\VV{k}{l}=\wedge^k\tzus\otimes\wedge^l\tuz\otimes\wedge^n\tuzs$ so that 
\[ \wedge^m\big(\tzu\oplus\tuzs\big)^*\otimes
\wedge^n\tuzs = 
\bigoplus_{k+l=m} \VV{k}{l} .\]
Since $A:=\tzu$ and $B:=\tuzspi$ are complex Lie subalgebroids of $\tzu\bowtie\tuzspi$, one has \[ \dee^{\nabla}_{\bowtie}\sections{\VV{k}{l}} 
\subset\sections{\VV{k+1}{l}\oplus\VV{k}{l+1}} .\] 
Composing $\dee^{\nabla}_{\bowtie}$ with the natural projections on each of the direct summands, we get the commutative diagram 
\[ \xymatrix{
&\sections{\VV{k}{l}} \ar[d]^{d_{\bowtie}^{\nabla}}
\ar[dl]_{\partial^{\nabla}_A} \ar[dr]^{(-1)^k\partial^{\nabla}_B} & \\
\sections{\VV{k+1}{l}} & \sections{\VV{k+1}{l}\oplus\VV{k}{l+1}} 
\ar[l] \ar[r] & 
\sections{\VV{k}{l+1}} ,} \]
where the operators $\partial^{\nabla}_A$ and $\partial^{\nabla}_B$ 
are given by
\beq{North_Carolina}\begin{aligned}
& \big(\partial^{\nabla}_A\alpha\big)(A_0,\dots,A_k,B_1,\dots,B_l) \\
= \; & \sum_{i=0}^k (-1)^i \Big(\nabla_{A_i}(\alpha(A_0,\dots,\widehat{A_i},\dots,A_k,
B_1,\dots,B_l)) \\
& - \sum_{j=1}^l \alpha(A_0,\dots,\widehat{A_i},\dots,A_k,
B_1,\dots, \pr_B\lb{A_i}{B_j},\dots,B_l)\Big) \\
& + \sum_{i<j} (-1)^{i+j} \alpha(\lie{A_i}{A_j},A_0,\dots,\widehat{A_i},
\dots,\widehat{A_j},\dots,A_k,B_1,\dots,B_l)
\end{aligned}\eeq
and
\beq{South_Carolina}\begin{aligned}
& \big(\partial^{\nabla}_B\alpha\big)(A_1,\dots,A_k,B_0,\dots,B_l) \\
= \; & \sum_{i=0}^l (-1)^i \Big(\nabla_{B_i}(\alpha(A_1,\dots,A_k,
B_0,\dots,\widehat{B_i},\dots,B_l)) \\
& - \sum_{j=1}^k \alpha(A_1,\dots,\pr_A\lb{B_i}{A_j},\dots,A_k,
B_0,\dots,\widehat{B_i},\dots,B_l)\Big) \\
& + \sum_{i<j} (-1)^{i+j} \alpha(A_1,\dots,A_k,
\lie{B_i}{B_j},B_0,\dots,\widehat{B_i},\dots,\widehat{B_j},\dots,B_l) ,
\end{aligned}\eeq
for all $\alpha\in\sections{\wedge^k A^*\otimes\wedge^l B^*}$, $A_0,\dots,A_k\in\sections{A}$ and $B_0,\dots,B_k\in\sections{B}$.
Here $\pr_B\lb{A_i}{B_j}$ denotes the $B$-component of $\lb{A_i}{B_j}\in A\bowtie B$ and $\pr_A\lb{B_i}{A_j}$ the $A$-component of $\lb{B_i}{A_j}$.

Since $d_{\bowtie}^{\nabla}=\partial_A^{\nabla}
+(-1)^k\partial_B^{\nabla}$, it follows from $(d_{\bowtie}^{\nabla})^2=0$ that 
$(\partial_A^{\nabla})^2=0$, $(\partial_B^{\nabla})^2=0$ and 
$\partial_A^{\nabla}\rond\partial_B^{\nabla}
=\partial_B^{\nabla}\rond\partial_A^{\nabla}$. 
Thus the complex~\eqref{Virginia} is the total complex of the double complex 
\begin{equation*}
\xymatrix{ \sections{\VV{k}{l}}
\ar[r]^{\partial^{\nabla}_B} \ar[d]_{\partial^{\nabla}_A} &
\sections{\VV{k}{l+1}} \ar[d]^{\partial^{\nabla}_A} \\
\sections{\VV{k+1}{l}} \ar[r]_{\partial^{\nabla}_B} &
\sections{\VV{k+1}{l+1}} }
\end{equation*}

Hence it follows that 
$H^*(\tzu\bowtie\tuzspi,\wedge^n\tuzs)$
is isomorphic to the total cohomology of the double complex 
\[ \xymatrix{ 
\sections{\wedge^i(\tzu)^*\otimes\wedge^j(\tuz)
\otimes\wedge^n(\tuz)^*} \ar[r]^{\del_{B}^{\nabla}} 
\ar[d]_{\del_{A}^{\nabla}} & 
\sections{\wedge^i(\tzu)^*\otimes\wedge^{j+1}(\tuz)
\otimes\wedge^n(\tuz)^*} \ar[d]^{\del_{A}^{\nabla}} \\ 
\sections{\wedge^{i+1}(\tzu)^*\otimes\wedge^j(\tuz)
\otimes\wedge^n(\tuz)^*} \ar[r]_{\del_{B}^{\nabla}} & 
\sections{\wedge^{i+1}(\tzu)^*\otimes\wedge^{j+1}(\tuz)
\otimes\wedge^n(\tuz)^*} 
} \] 
By $\tau$ we denote the natural contraction map
\beq{Georgia} \tau: \sections{(\wedge^i(\tzu)^*\otimes\wedge^j(\tuz) )^*
\otimes\wedge^n(\tuz)^*}\to \Omega^{n-j, i} ,\eeq
which is an isomorphism of $\FF$-modules.

Take a local holomorphic chart $(U;z_1,\dots,z_n)$ of $X$, and set 
\begin{gather*}
b=\pz{j_1}\wedge\cdots\wedge\pz{j_l} ;\\ 
\omega=\dz{1}\wedge\cdots\wedge\dz{n} 
.\end{gather*} 
Because of~\eqref{Iowa}, we have
\beq{Kentucky} 
\tau d_{\pi}^{\nabla} (b\otimes \omega)=(-1)^{l+1} \del_{\pi} \tau (b\otimes \omega) 
.\eeq

\begin{lem}\label{Detroit}
For all $f\in\FF$, $b\in\hXX^k$ and $\mu\in\hOO^l$, we have: 
\begin{gather*} 
d_{\pi}(fb)=-(\pi\diese\del f)\wedge b+f(d_{\pi}b) , \\ 
\delpi(f\mu)=(\pi\diese\del f)\brc\mu+f(\delpi\mu) .
\end{gather*}
\end{lem}

As a consequence, we have
\begin{prop}
\begin{gather} 
\tau\rond\partialA=\delbar\rond\tau \label{Delaware} \\ 
\tau\rond\partialB=(-1)^{k+l+1}\delpi\rond\tau \label{Maine} 
\end{gather}
\end{prop}

\begin{proof}
The first relation~\eqref{Delaware} is a simple consequence of the definition~\eqref{North_Carolina} of $\partialA$, while the second~\eqref{Maine} follows from~\eqref{South_Carolina}, \eqref{Kentucky} and Lemma~\ref{Detroit}. 
\end{proof}

Now the conclusion of the theorem follows immediately.
\end{proof}

\section{Evens-Lu-Weinstein duality}

We recall a remarkable duality construction due to Evens, Lu \& Weinstein~\cite{ELW}.

Consider a compact complex (and therefore orientable) manifold $X$ with $\dim_{\CC}X=n$, a complex Lie algebroid $B$ over $X$ with $\rk_{\CC}B=r$ and a module $E$ over $B$. The complex dual $E^*$ is also a module over $B$. 
We will use the symbol $\nabla$ to denote the representations of $B$ on both $E$ and $E^*$. 

The complex vector bundle $Q_B=\wedge^r B\otimes \wedge^{2n}\tcxs$ 
is a module over  
the complex Lie algebroid $B$ with representation
 $D:\sections{Q_B}\to\sections{B^*\otimes Q_B}$~\cite{ELW} given by 
\[ D_b(X\otimes\mu)=\lb{b}{X}\otimes\mu+X\otimes\ld{\rho(b)}\mu ,\]
for all $b\in\sections{B}$, $X\in\sections{\wedge^r B}$ and $\mu\in\sections{\wedge^{2n}\tcxs}$. 

By $H^*(B,E)$ and $H^*(B,E^*\otimes Q_B)$, we denote 
the  Lie algebroid cohomology of $B$ with coefficients
in $E$ and $E^*\otimes Q_B$, respectively.
We use the  notation $\dnb$ to denote their coboundary differential
operators in both cases.
Let $\Xi$ be the isomorphism of vector bundles:
\begin{equation*}
\Xi:\wedge^r B^*\otimes(\wedge^r B\otimes\wedge^{2n}\tcxs) \to \wedge^{2n}\tcxs :
\xi\otimes(X\otimes\mu) \mapsto (\xi\brc X)\mu
.\end{equation*}

The following lemma can be verified by a direct computation.

\begin{lem}
We have 
\[ \Xi\rond\dnb\big(\xi\otimes(X\otimes\mu)\big) 
=(-1)^{r-1} d\big(\rho(\xi\brc X)\brc\mu\big) ,\]
for any $\xi\otimes(X\otimes\mu  )\in \sections{\wedge^{r-1} B^*\otimes Q_B}$
\end{lem}

Consider the bilinear map 
\[ \prepairing{\cdot}{\cdot}:\sections{\wedge^k B^* \otimes E}\otimes 
\sections{\wedge^{r-k}B^*\otimes E^*\otimes Q_B}
\to\sections{\wedge^{2n}\tcxs} \]
defined by 
\[ \prepairing{\xi_1\otimes e}{\xi_2\otimes\epsilon\otimes(X\otimes\mu)} 
= \epsilon(e) \cdot (\xi_1\wedge\xi_2)(X) \cdot \mu .\]

\begin{lem}
If $\xi_1\otimes e\in\sections{\wedge^{k-1}B^*\otimes E}$ and 
$\xi_2\otimes\epsilon\otimes(X\otimes\mu)\in\sections{\wedge^{r-k}B^*\otimes E^*\otimes Q_B}$, then 
\begin{multline*} 
\prepairing{\dnb(\xi_1\otimes e)}{\xi_2\otimes\epsilon\otimes(X\otimes\mu)} 
+ (-1)^{r-1}\prepairing{\xi_1\otimes e}{\dnb\big(\xi_2\otimes\epsilon\otimes(X\otimes\mu)\big)} \\ 
= \Xi\rond\dnb\big(\epsilon(e)\cdot\xi\otimes (X\otimes\mu)\big)
= (-1)^{r-1} d\big(\epsilon(e)\cdot\rho(\xi\brc X)\brc \mu\big)
.\end{multline*}
\end{lem}

Therefore, by Stokes' theorem, the pairing 
\[ \pairing{\alpha }{\beta} 
= \int_X \prepairing{\alpha}{\beta} \] 
where $\alpha\in \sections{\wedge^k B^*\otimes E}$ and 
$\beta\in \sections{\wedge^{r-k} B^*\otimes E^*\otimes Q_B}$
satisfies 
\[ \pairing{\dnb(\alpha)}{\beta}
+ (-1)^{r-1} \pairing{\alpha}{\dnb(\beta)} = 0 \] (where $\alpha\in\sections{\wedge^{k-1} B^*\otimes E}$ and 
$\beta\in \sections{\wedge^{r-k} B^*\otimes E^*\otimes Q_B}$)  
and thus induces a pairing at the cohomology level~\cite{ELW}: 
\begin{equation}
\label{Wyoming}
 \pairing{\cdot}{\cdot}:H^k(B,E)\otimes H^{r-k}(B,E^*\otimes Q_B) 
\to \CC .
\end{equation}

The following is due to Block~\cite{Block}.

\begin{prop}
If  $B$ is an elliptic Lie algebroid, the cohomology pairing~\eqref{Wyoming} is perfect.
\end{prop}

Given a holomorphic Poisson manifold $(X,\pi)$, we can take 
$B=\tzu\bowtie\tuzspi$. Then $Q_B^{\half}=\wedge^n \tuzs$ and, taking $E=Q_B^{\half}$, we have $E=\wedge^n\tuzs$ and $E^*\otimes Q_B=\wedge^n\tuzs$. 

In this particular case, we get the cohomology pairing 
\[ \pairing{\cdot}{\cdot}: 
H^k\big(\tzu\bowtie\tuzspi,\wedge^n\tuzs\big) \otimes 
H^{2n-k}\big(\tzu\bowtie\tuzspi,\wedge^n\tuzs\big) 
\to \CC ,\] 

If we identify the cochain group $\bigoplus_{k,l}\VV{k}{l}$ with $\bigoplus_{p,q}\hOO^{p,q}$ via the contraction map $\tau$ (see Equation~\eqref{Georgia}), then a straightforward (though lengthy) computation 
shows that, on the cochain level, the above cohomology pairing is given by 
\[ \hOO^{i,j}\otimes\hOO^{k,l}\to\CC: \zeta\otimes\eta\mapsto\int_X 
(\zeta\wedge\eta)^{top} .\] 

We have proved the following theorem.

\begin{thm}
Let $(X,\pi)$ be  a compact holomorphic Poisson manifold.
The pairing 
\[ \pairing{\cdot}{\cdot}:H_{2n-k}(X,\pi)\otimes H_k(X,\pi)\to\CC:
[\zeta]\otimes[\eta]\mapsto\int_X (\zeta\wedge\eta)^{top} \] 
(where $\zeta,\eta\in\oplus_{k,l}\hOO^{k,l}$) 
is nondegenerate.
\end{thm}

\begin{rmk}
When $X$ is a compact complex manifold considered
as a zero Poisson manifold, 
then $H_k(X,\pi) \cong \oplus_{j-i=n-k} H^{i, j} (X)$.
The above theorem easily follows from Serre duality.
\end{rmk}

\section{Examples}

The purpose of this section is the computation of the Koszul-Brylinski Poisson homology
of all Poisson structures with which $\cposq$ can be endowed. 

From now on, $X$ will denote the complex manifold $\cposq$. 
Since $X$ is 2-dimensional, all holomorphic bivector fields on it are automatically Poisson tensors.
Thus the Poisson tensors on $X$ form the complex vector space $H^0(X,\wedge^2\tx)$, which is know to be 9-dimensional.
Here is a more explicit description of $H^0(X,\wedge^2\tx)$.

\begin{prop}
\label{South_Bend}
\cite{Duistermaat}
Let $P^{2,2}$ denote the $9$-dimensional vector space of all
bihomogeneous polynomials on $\CC^2\times \CC^2$ of bidegree $(2,2)$.
Given any $p\in P^{2,2}$, there exists a unique holomorphic
bivector field $\pi_p$ on $X=\cposq$ such that, 
in an affine chart $(z_1,z_2)\mapsto([1:z_1],[1:z_2])$ of $\cposq$, we have
\[ \pi_p=q(z_1,z_2)\;\partial_{z_2}\wedge\partial_{z_1} ,\]
where $q(z_1,z_2)=p\big((1,z_1),(1,z_2)\big)$. The map 
\[ p\in P^{2,2}\longmapsto \pi_p\in H^0 (X, \wedge^2 T_X) \]
is an isomorphism of complex vector spaces.
As a consequence, the space of all holomorphic Poisson
 bivector fields is a $9$-dimensional vector space over $\CC$.
\end{prop}

\begin{thm}
For any holomorphic Poisson bivector field $\pi$ on $X=\cposq$, we have 
\[ H_0(X,\pi)=0, \; H_1(X,\pi)=0, \; H_2(X,\pi)\cong\CC^4, \; H_3(X,\pi)=0, \; H_4 (X,\pi)=0 .\]
\end{thm}

\begin{proof}
Let us first assume that $\pi=0$.
In this case, $H_k(X,\pi)=\oplus_{j-i=2-k}H^{i,j}(X)$ 
and we obtain
\begin{align*} 
H_0(X,\pi)&=H^{0,2}(X)=0, &
H_1(X,\pi)&=H^{1,2}(X)\oplus H^{0,1}(X)=0, \\ 
H_3(X,\pi)&=H^{2,0}(X)=0, & 
H_4(X,\pi)&=H^{2,1}(X)\oplus H^{1,0}(X)=0, 
\end{align*}
and  
\[ H_2(X,\pi)=H^{0,0}(X)\oplus H^{1,1}(X)\oplus H^{2,2}(X)\cong\CC^4 .\]

Now let us assume that $\pi\neq 0$.
By definition, $H_0(X,\pi)$ consists of those $\alpha \in \Omega^{2,0}_X$ such that
$\delbar \alpha =0$ and $\partial_\pi\alpha=0$.
The first condition means that $\alpha$ is a holomorphic
2-form on $X$. Since $H^0(X,\hOO^2)=H^{0,2}(X)
=H^{0,2}(\cposq)=0$, it follows that
$H_0(X,\pi)=0$.

We now proceed to compute $H_1(X,\pi)$. Assume that $\theta+\omega\in\Omega^{1,0}\oplus\Omega^{2,1}$ is a Koszul-Brylinski 1-cycle.
That is, $\partial_\pi\theta+\delbar \theta+\partial_\pi\omega+\delbar \omega=0$.
Hence it follows that $\partial_\pi \theta=0$, $\delbar \theta +\partial_\pi \omega=0$ and  $\delbar \omega=0$. 
Since $H^{1,2}(\cposq)=0$, there exists $\beta\in\Omega^{2,0}_X$ with $\omega=\delbar\beta$. 
On the other hand, from $\partial_\pi\theta=0$, it follows that $i_{\pi}\partial \theta=0$ since $\partial_\pi\theta=[\partial,i_\pi]\theta =-i_{\pi}\partial\theta$. Therefore, $\partial\theta$ vanishes at those points where 
$\pi$ does not vanish. By Proposition~\ref{South_Bend},
$\pi$ is nonzero on a dense subset of $\cposq$.
Thus, we have $\partial\theta=0$, which implies that
$\theta=\partial\alpha$ for some $\alpha\in\Omega^{0,0}_X$
since $H^{1,0}(\cposq)=0$.
It follows that 
\[ 0=\delbar \theta +\partial_\pi \omega=
\delbar \partial \alpha+\partial_\pi \delbar \beta
=\delbar  (-\partial \alpha-\partial_\pi \beta) .\]
Since $H^{0, 1}(\cposq)=0$, we have $\partial\alpha+\partial_\pi\beta=0$. 
Thus $\theta+\omega=(\delbar-\partial_\pi)\beta$ from which we conclude that $H_1(X,\pi)=0$.

By Evens-Lu-Weinstein duality, we have 
\[ H_3 (X, \pi)\cong  H_1 (X, \pi)=0 \qquad \text{and} \qquad H_4 (X, \pi)\cong H_0(X, \pi)=0 .\]
Moreover, according to Theorem~\ref{Silver_Springs},
\[ \EulerKB{X,\pi}=\Euler{X}=\Euler{\CP{1}}+\Euler{\CP{1}}=4 .\] 
Thus we have $H_2(X,\pi)\cong\CC^4$.
This concludes the proof.
\end{proof}

\end{document}